\documentclass[11pt,twoside]{article}

\usepackage{amsfonts}
\usepackage{amsbsy}
\usepackage{amssymb}
\usepackage{amsmath}
\usepackage{amsthm}
\usepackage{indentfirst}
\usepackage[pdftex]{graphicx}
\usepackage{subfigure}
\usepackage{supertabular}
\usepackage{algorithm}
\usepackage{enumitem}
\usepackage{epstopdf}
\usepackage{listings}
\usepackage{color}
\usepackage{courier}
\usepackage{algorithm}
\usepackage{algpseudocode}
\usepackage{verbatim} 
\usepackage[mathcal]{eucal}
\usepackage{multicol}
\usepackage{float}

\newtheorem{definition}{Definition}
\newtheorem{theorem}{Theorem}
\newtheorem{lemma}{Lemma}

\newtheorem{remark}{Remark}
\newtheorem{proposition}{Proposition}
\begin{document}

\begin{center}
{\Large The weak $n$-inner product space} 
\end{center}
\medskip
\begin{center}

Nicu\c{s}or Minculete        and
        Radu P\u alt\u anea 
\end{center}

\begin{abstract}
In this article we study a generalization of the $n$-inner product which we name weak $n$-inner product. As particular case we consider the $n$-iterated $2$-inner product and we give its representation in terms of the standard $k$-inner products, $k\le n$, using the Dodgson's identity for determinants. Finally, we present several applications, including a brief characterization of a linear regression model for the random variables in discrete case and  a generalization of the Chebyshev functional using the $n$-iterated $2$-inner product.\\


\end{abstract}

\section{Introduction}\label{sec:1}

\indent The concept of linear 2-normed spaces and 2-metric spaces has been investigated by G\" ahler \cite{8}. In \cite{3} and \cite{4}, Diminnie, G\" ahler and White studied the 2-inner product spaces.\\
\indent A classification of results related to the theory of 2-inner product spaces can be found in book \cite{1}. Here, several properties of 2-inner product spaces are given. In \cite{6} Dragomir et al. show the corresponding version of Boas-Bellman inequality in 2-inner product spaces. Others properties of a 2-inner product space can be found in \cite{2}.\\
\indent Misiak \cite{11} generalizes this concept of a 2-inner product space, in 1989, in the following way: let $n$ be a nonnegative integer $(n \geq 2)$ and $X$ be a vector space of dimension $\dim X=d\geq n$   ($d$ may be infinite) over the field of real numbers $\mathbb{R}$. An $\mathbb{R}$-valued function $\langle\cdot,\cdot\mid\cdot,...,\cdot\rangle$ on $X^{n+1}$ satisfying the following properties:\\
\indent I1) $\langle v_{1},v_{1}| v_{2},...,v_{n}\rangle\geq 0; \langle v_{1},v_{1}|v_{2},...,v_{n}\rangle= 0$ if and only if $v_{1},v_{2},...,v_{n}$ are linearly dependent;\\
\indent I2) $\langle v_1,v_1| v_{2},...,v_{n}\rangle=\langle v_{i_1},v_{i_1}| v_{i_2},v_{i_3},...,v_{i_n}\rangle$, for every permutation\linebreak $(i_1,i_2,...,i_n)$ of $(1,2,...,n)$;\\
\indent I3) $\langle v,w| v_{2},...,v_{n}\rangle=\langle w,v| v_{2},...,v_{n}\rangle$;\\
\indent I4) $\langle\alpha v,w| v_{2},...,v_{n}\rangle=\alpha\langle v,w| v_{2},...,v_{n}\rangle$,  for every scalar $\alpha \in \mathbb{R}$.\\
\indent I5) $\langle v+v',w| v_{2},...,v_{n}\rangle=\langle v,w| v_{2},...,v_{n}\rangle+\langle v',w| v_{2},...,v_{n}\rangle$;\\
is called an \textit{ n-inner product} on $X$, and the pair $(X,\langle\cdot,\cdot|\cdot,...,\cdot\rangle)$ is called an \textit{n-inner product space} or \textit{n-pre-Hilbert space}.\\
\indent It is easy to see that the $n$-inner product is a linear function of its two first arguments. Several results related to the theory of the $n$-inner product spaces can be found in \cite{9}, \cite{12}: $\langle v,w|\alpha v_{2},...,v_{n}\rangle=\alpha^{2}\langle v,w| v_{2},...,v_{n}\rangle$, 
for every real number $\alpha$ and for $v,w,v_{2},..,v_{n} \in X$;
 $\langle v,w|v_{2}+v_{2}',v_{3},...,x_{n}\rangle-\langle v,w|v_{2}-v_{2}',v_{3},...,v_{n}\rangle=\langle v_{2},v_{2}'|v+w,v_{3},...,v_{n}\rangle-\langle v_{2},v_{2}'|v-w, v_{3},...,v_{n}\rangle$, for all $v,w,v_{2},v_{3},...,v_{n}, v_{2}' \in X$ and an extension of the Cauchy-Schwarz inequality to arbitrary \textit{n}:\\ 
\begin{equation}\label{1}
|\langle v,w|v_{2},...,v_{n}\rangle|\leq \sqrt{\langle v,v|v_{2},...,v_{n}\rangle}\sqrt{\langle w,w|v_{2},...,v_{n}\rangle},
\end{equation}
for all $v,w,v_{2},...,v_{n} \in X$. The equality holds in (1) if and only if $v,w,v_{2},...,v_{n}$  are linearly dependent.\\
\indent Other consequences from the above properties can be inferred very easily: 
$$\langle 0,w|v_{2},...,v_{n}\rangle=\langle v,0|v_{2},...,v_{n}\rangle=\langle v,w|0,...,v_{n}\rangle=0,$$
$$ \langle v_{2},w|v_{2},...,v_{n}\rangle=\langle v,v_{2}|v_{2},...,v_{n}\rangle=0,$$
for all $v,w,v_{2},...,v_{n}\in X$.\\
\indent Let $(X,\langle\cdot,\cdot|\cdot,...,\cdot\rangle)$ be an $n$-inner product space, $n\geq2$. We can define a function $\|\cdot,...,\cdot\|$ on $X\times X\times ...\times X=X^{n}$ by
$$\|v|v_{2},...,v_{n}\|:=\sqrt{\langle v,v|v_{2},...,v_{n}\rangle},$$
for all $v,v_{2},...,v_{n} \in X$, which in \cite{11}  is shown that satisfies the following conditions: \\
\indent N1) $\|v|v_{2},...,v_{n}\|\geq0$ and $\|v|v_{2},...,v_{n}\|=0$ if and only if $v,v_{2},...,v_{n}$ are linearly dependent;\\
\indent N2) $\|v|v_{2},...,v_{n}\|$ is invariant under permutation;\\
\indent N3) $\|\alpha v|v_{2},...,v_{n}\|=|\alpha|\|v|v_{2},...,v_{n}\|$, for any scalar $\alpha \in \mathbb{R}.$\\
\indent N4) $\|v+w|v_{2},...,v_{n}\|\leq \|v|v_{2},...,v_{n}\|+\|w|v_{2},...,v_{n}\|$,\\
 for all $v,w,v_{2},...,v_{n}\in X.$
 
 A function $\| \cdot|\cdot,...,\cdot\|$  defined on $X^{n}$ and satisfying the above conditions is called an \textit{n-norm} on X and $(X,\| \cdot|\cdot,...,\cdot\|  )$ is called a linear \textit{n-normed space}.\\
\indent It is easy to see that if $(X,\langle \cdot,\cdot|\cdot,...,\cdot\rangle)$  is an \textit{n}-inner product space over the field of real numbers $ \mathbb{R}$, then $(X,\| \cdot|\cdot,...,\cdot\|  )$ is a linear \textit{n}-normed space and the \textit{n}-norm $\| \cdot|\cdot,...,\cdot\|  $  is generated by  the \textit{n}-inner product $\langle\cdot,\cdot|\cdot,...,\cdot\rangle$.\\
\indent Furthermore, we have the parallelogram law \cite{1},\\
\begin{equation}\label{2}
\|v+w|v_{2},...,v_{n}\|^{2}+\|v-w|v_{2},...,v_{n}\|^{2}=2\|v|v_{2},...,v_{n}\|^{2}+2\|w|v_{2},...,v_{n}\|^{2},
\end{equation}
for all $v,w,v_{2},..,v_{n}\in X$ and the polarization identity (see e.g. \cite{1} and \cite{2}),\\
\begin{equation}\label{3}
\| v+w|v_{2},...,v_{n}\|^{2}-\| v-w|v_{2},...,v_{n}\|^{2}=4\langle v,w|v_{2},...,v_{n}\rangle,
\end{equation}
for all $v,w,v_{2},...,v_{n} \in X.$\\
The {\it standard $n$-inner product} on an inner product space $X=(X,\langle\cdot.\cdot\rangle)$ is given by:
\begin{equation}\label{4}
\langle v,w|v_{2},...,v_{n}\rangle:= \begin{vmatrix}
\langle v,w\rangle & \langle v,v_{2}\rangle &... & \langle v,v_{n}\rangle \\ \langle v_{2},w\rangle & \langle v_{2},v_{2} \rangle &... & \langle v_{2},v_{n}\rangle \\
\vdots & \vdots &\vdots & \vdots \\
 \langle v_{n},w \rangle & \langle v_{n},v_{2} \rangle &... & \langle v_{n},v_{n}\rangle
\end{vmatrix},
\end{equation}
which generates $n$-norm $\| v|v_{2},...,v_{n}\|:=\sqrt{\langle v,v|v_{2},...,v_{n}\rangle},$ representing the volume of the $n$-dimensional parallelepiped spanned by $v,v_{2},...,v_{n}.$
\par Various type of applications of $n$-inner products and $n$-norms can be found in recent papers \cite{BC}, \cite{G18}, \cite{HT}, \cite{Ka}, \cite{MSA}, \cite{MS}, \cite{RJ16}, \cite{RMA}.

\begin{remark} \label{R1.1}
{\rm The standard $n$-inner product satisfies also the following additional condition:\\
\indent I6) If $v,v_2,\ldots,v_n$ are linearly dependent, then $\langle v,w|v_2,\ldots,v_n\rangle=0$,\\
for $v,w,v_2,\ldots,v_n\in X$.}
\end{remark}

\indent The motivation of this article is to study another type of $n$-inner product built based on the properties of the $n$-inner product, except property I2. We will define the weak $n$-inner product and the  $n$-iterated $2$-inner product and we will give its representation in terms of the  standard  $k$-inner products, $k\le n$, using the Dodgson's identity for determinants. We also present a brief characterization of a linear regression model for the random variables in discrete case. Finally, we generalize the Chebyshev functional using the $n$-iterated $2$-inner product.

\section{The weak $n$-inner product}\label{sec:2}

Let $X$ be a real vector space. 

\begin{definition}\label{D1}
\par An $\mathbb{R}$-valued function $(\cdot,\cdot\mid\cdot,...,\cdot)$ on $X^{n+1}$, $n\ge 2$, satisfying the following properties:
\par{\rm P1)} {\it Positivity}: $(x,x| x_{n},...,x_{2})\geq 0\;\mbox{and}\; (x,x|x_{n},...,x_{2})= 0$ if and only if $x, x_{2},x_{3},...,x_{n}$ are linearly dependent;
\par{\rm P2)} {\it Interchangeability}: $(x,x| x_{n},...,x_{2})=(x_{n},x_{n}| x,x_{n-1},...,x_{2})$;
\par{\rm P3)} {\it Symmetry}: $(x,y| x_{n},...,x_{2})=(y,x| x_{n},...,x_{2})$;
\par{\rm P4)} {\it Homogeneity}: $(\alpha x,y| x_{n},...,x_{2})=\alpha(x,y| x_{n},...,x_{2})$,  for every scalar $\alpha \in \mathbb{R}$.
\par{\rm P5)} {\it Additivity}: $(x+x',y| x_{n},...,x_{2})=(x,y| x_{n},...,x_{2})+(x',y| x_{n},...,x_{2})$;\\
is called a {\bf weak $n$-inner product} on $X$, and the pair $(X,(\cdot,\cdot|\cdot,...,\cdot))$ is called a {\bf weak $n$-inner product space} or {\bf weak n-pre-Hilbert space}.
\end{definition}

\begin{remark} {\rm It is easy to see that: 
$$(0,y|x_{n},...,x_{2})=(x,0|x_{n},...,x_{2})=(x,y|0,...,x_{2})=0.$$}
\end{remark}

\begin{remark} {\rm Obviously an $n$-inner product is a weak $n$-inner product, so an $n$-inner product space is a weak $n$-inner product space, but the reciprocal is not true. This fact will be shown in Remark \ref{R2.4}.
\par For $n=2$ a weak $n$-inner product is also an $n$-inner product. For $n\ge 3$ a weak $n$-inner product can be build, for instance, by formula
$$(x,y|x_{n},\ldots,x_2)=\Theta(x,y|x_{n}\ldots,x_2)\cdot \Psi(x_{n-1},\ldots,x_{2}),$$
where $\Theta(x,y|x_{n}\ldots,x_2)$ is a $n$-inner product} and
 $\Psi:X^{n-2}\to\mathbb{R}$ is a function with properties  $\Psi(x_{n-1},\ldots,x_{2})\ge 0.\;\forall (x_{n-1},\ldots,x_{2})$ and $\Psi(x_{n-1},\ldots,x_{2})=0$ iff $x_{n-1},\ldots,x_{2}$ are linearly dependent (in the case $n=3$, this means $x_2=0$).
\end{remark}

In the next lemma we generalize a property that exists in the case of $2$-inner products. The method of the proof is based on the method used in \cite{2}.

\begin{lemma}\label{L1a} Let $,x_2,\ldots,x_n,x,y\in X$. If $x_2,\ldots,x_n,x$ are linearly dependent, then
\begin{equation}\label{eP6}
(x,y|x_n,\ldots,x_2)=0.
\end{equation}
\end{lemma}
\par{\it Proof} We consider two cases.
\par \underline{Case 1.} $x_2,\ldots,x_n,y$ {\it are linearly independent}. Consider the vector
$$u=(y,y|x_n,\ldots,x_2)x-(x,y|x_n,\ldots,x_2)y.$$
Then from P1) we hve $(u,u|x_n,\ldots,x_2)\ge 0$. This inequality is equivalent to
$$(y,y|x_n,\ldots,x_2)[(x,x|x_n,\ldots,x_2)(y,y|x_n,\ldots,x_2)-(x,y|x_n,\ldots,x_2)^2]\ge 0.$$
Since $x_2,\ldots,x_n,x$ are linearly dependent, from P1) we obtain $(x,x|x_n,\ldots,x_2)=0$ and hence
$$-(y,y|x_n,\ldots,x_2)(x,y|x_n,\ldots,x_2)^2\ge 0.$$
Since  $x_2,\ldots,x_n,y$ are linearly independent it follows that $(y,y|x_n,\ldots,x_2)>0$. Consequently one obtains (\ref{eP6}).
\par \underline{Case 2.} $x_2,\ldots,x_n,y$ {\it are linearly dependent}. Then also $x_2,\ldots,x_n,x+y$ are are linearly dependent. We have
$$(x,y|x_n,\ldots,x_2)=\frac12[(x+y,x+y|x_n,\ldots,x_2)-(x,x|x_n,\ldots,x_2)-(y,y|x_n,\ldots,x_2)].$$
Because $(x,x|x_n,\ldots,x_2)=0$, $(y,y|x_n,\ldots,x_2)=0$, $(x+y,x+y|x_n,\ldots,x_2)=0$ relation (\ref{eP6}) follows.
\framebox

\begin{theorem}\label{T0} 
Suppose that $(X,( \cdot,\cdot|\cdot,...,\cdot))$  is a weak \textit{n}-inner product space over the field of real numbers $ \mathbb{R}$. Let $x_2,\ldots, x_n\in X$, $n\ge 2$ be fixed. Denote $Y=span\{x_2,\ldots,x_n\}$. Define the quotient space $X/Y=\{\hat{x}|\;x\in X\}$, where $\hat{x}=\{u\in X|u-x\in Y\}$, $x\in X$. Then function $\psi: (X/Y)^2\to\mathbb{R}$, $\psi(\hat{x},\hat{y}):= (x,y|x_{n},...,x_{2})$, 
 $\hat{x},\hat{y}\in X/Y$ is well defined and is a semi-inner product on $X/Y$.  Moreover, if $x_2,\ldots,x_n$ are linearly independent, then $\psi$ is an inner product.
\end{theorem}
\par{\it Proof} Let $x,x',y,y'\in X$, such that $x'-x\in Y$ and $y'-y\in Y$. Using Lemma \ref{L1a} we get $\psi(\widehat{x'},\hat{y'})=(x',y'|x_{n},...,x_{2})=
(x,y|x_{n},...,x_{2})+(x'-x,y|x_{n},...,x_{2})+(x,y'-y|x_{n},...,x_{2})+(x'-x,y'-y|x_{n},...,x_{2})=(x,y|x_{n},...,x_{2})=\psi(\widehat{x},\hat{y})$.
This means that $\psi$ is well defined.
\par From {\rm P1)} we have $\psi(\hat{x},\hat{x})=(x,x|x_{n},...,x_{2})\ge0$.  Moreover, if $\psi(\hat{x},\hat{x})=0$, then $(x,x|x_2,\ldots,x_n)=0$, which implies that $x,x_2,\ldots,x_n$ are linearly dependent. If $x_2,\ldots, x_n$ are linearly independent it follows that $x\in Y$.  Then $\hat{x}=\hat{0}$. 
\par The other properties of the inner product follow in a simple manner from conditions {\rm P3)}, {\rm P4)} and {\rm P5)}.
\framebox

\begin{theorem}\label{T1} {\bf ( Schwarz type inequality)}
Let $(X,(\cdot,\cdot|\cdot,...,\cdot))$ be a weak $n$-inner product space. For any $x,y,x_{2},...,x_{n} \in X$ we have
\begin{equation}\label{eSchwarz}
|(x,y|x_{n},...,x_{2})|\leq \sqrt{(x,x|x_{n},...,x_{2})}\sqrt{(y,y|x_{n},...,x_{2})}.
\end{equation}
 In the case when $x_2,\ldots,x_n$ are linearly independent, then the equality holds in (\ref{eSchwarz}) if and only if there exist $\mu\in \mathbb{R}_+$ and $u\in Y:=span\{x_2,\ldots,x_n\}$ such that $y=\mu x+u$.
\end{theorem}
\par{\it Proof} 
 By taking into account Theorem \ref{T0} and the notations given there, we have
\begin{eqnarray*}
&&|(x,y|x_{n},...,x_{2})|\\
&&\qquad=|\Psi(\hat{x},\hat{y})|\le \sqrt{\Psi(\hat{x},\hat{x})}\sqrt{\Psi(\hat{y},\hat{y})}=\sqrt{(x,x|x_{n},...,x_{2})}\sqrt{(y,y|x_{n},...,x_{2})}.
\end{eqnarray*}
If $x_2,\ldots,x_n$ are linearly independent, then the equality holds in (\ref{eSchwarz}) iff there is $\mu\ge 0$, such that $\hat{y}=\mu\hat{x}$, i.e. exists $u\in Y$ for which $y=\mu x+u$. 
\framebox.

\begin{definition}\label{D2}
Let $(X,(\cdot,\cdot|\cdot,...,\cdot))$ be a weak $n$-inner product space, $n\geq2$. We can define a function $\|\cdot,...,\cdot\|$ on $X\times X\times ...\times X=X^{n}$ by
\begin{equation}\label{eD2}
\|x|x_{n},...,x_{2}\|:=\sqrt{(x,x|x_{n},...,x_{2})},\; \mbox{for all}\; x,x_{2},...,x_{n} \in X.
\end{equation} 
\end{definition}

\begin{proposition}\label{P1}  
If  $(X,(\cdot,\cdot|\cdot,...,\cdot))$ is a weak $n$-inner product space, then function $\|\cdot,...,\cdot\|$ defined in (\ref{eD2}) 
satisfies the following conditions:
\begin{itemize}
\item[{\rm C1)}] $\|x|x_{n},...,x_{2}\|\geq0$ and $\|x|x_{n},...,x_{2}\|=0$ if and only if $x, x_{2},...,x_{n}$ are linearly dependent;
\item[{\rm C2)}]  $\|x|x_{n},x_{n-1},...,x_{2}\|=\|x_n|x,x_{n-1},...,x_{2}\|$;
\item[{\rm C3)}]  $\|\alpha x|x_{n},...,x_{2}\|=|\alpha|\|x|x_{n},...,x_{2}\|$, for any scalar $\alpha \in \mathbb{R};$
\item[{\rm C4)}]  $\|x+y|x_{n},...,x_{2}\|\leq \|x|x_{n},...,x_{2}\|+\|y|x_{n},...,x_{2}\|$
\end{itemize}
for all $x,y,x_{2},...,x_{n}\in X.$
\end{proposition}
\par{\it Proof} Conditions C1)-C4) follow immediately from conditions P1)-P5) and Definition \ref{D2}. \framebox

\begin{definition}\label{D3}
Let $X$ be a real vector space. A real function $\| \cdot|\cdot,...,\cdot\|$  defined on $X^{n}$ and satisfying conditions C1)-C4) is called a {\bf weak $n$-norm} on $X$ and $(X,\| \cdot|\cdot,...,\cdot\|  )$ is called a linear \textit{weak $n$-normed space}.
\end{definition}

\par It follows that if $(X,( \cdot,\cdot|\cdot,...,\cdot))$  is a weak \textit{n}-inner product space over the field of real numbers $ \mathbb{R}$, then $(X,\| \cdot|\cdot,...,\cdot\|  )$ is a linear weak \textit{n}-normed space and the weak \textit{n}-norm $\| \cdot|\cdot,...,\cdot\|  $  is generated by the weak \textit{n}-inner product $(\cdot,\cdot|\cdot,...,\cdot)$.

\begin{theorem}\label{T2} 
In conditions of Theorem \ref{T0}, function $\varphi:X/Y\to\mathbb{R}_+$, $\varphi(\hat{x}):= \|x|x_{n},...,x_{2}\|$, $\hat{x}\in X/Y$ is well defined and is a semi-norm on $X/Y$. Moreover, if $x_2,\ldots,x_n$ are linearly independent, then $\varphi$ is a norm.
\end{theorem}
\par{\it Proof} It follows immediately from Theorem \ref{T0}, since $\varphi(\hat{x})=\sqrt{\psi(\hat{x},\hat{x})}$, $\hat{x}\in X/Y$, where function $\psi$ was defined in this theorem. \framebox
\par\
\par In an inner product space, a special weak $n$-inner product  can be defined by recurrence starting from the $2$-inner product. Recall that the $2$-inner product was studied in \cite{1}, \cite{2}.

\begin{definition}\label{D3}
Let $(X,\langle \cdot,\cdot\rangle)$ be a real pre-Hilbert space. The $n$-{\bf iterated $2$-inner product}, or {\bf standard weak $n$-inner product} $(\cdot,\cdot|\cdot,\ldots,\cdot)_*:X^{n+1}\to\mathbb{R}$ is defined for $n\ge 2$ as follows. For $n=2$, $(\cdot,\cdot|\cdot)_*$ coincides with the standard $2$-inner product, i.e.
\begin{equation}\label{eD3a}
(x,y|z)_*:=\langle x,y|z\rangle=\begin{vmatrix} \langle x,y \rangle & \langle x,z \rangle \\ \langle z,y \rangle &\langle z,z \rangle \end{vmatrix}=\langle x,y\rangle\langle z,z \rangle-\langle x,z\rangle \langle z,y\rangle,\quad x,y,z \in X.
\end{equation}
Then, if $n\ge 3$ and $x,y,x_{2},...,x_{n}\in X$, define: 
\begin{equation}\label{eD3b}
(x,y|x_{n},...,x_{2})_*:=\begin{vmatrix}
(x,y|x_{n-1},...,x_{2})_* & (x,x_{n}|x_{n-1},...,x_{2})_*\\ (x_{n},y|x_{n-1},...,x_{2})_* &(x_{n},x_{n}|x_{n-1},...,x_{2})_*
\end{vmatrix}. 
\end{equation}
\end{definition}

\begin{theorem}\label{T3} 
If $(X,\langle \cdot,\rangle)$ is a real pre-Hilbert space, then for any $n\ge 2$ function  $(\cdot,\cdot|\cdot,\ldots,\cdot)_*:X^{n+1}\to\mathbb{R}$ given in Definition \ref{D3} is a weak $n$-inner product. 
\end{theorem}

\par{\it Proof} Consider proposition $S(n)$: the $n$-iterated $2$-inner product satisfies conditions P1)-P6). We prove this proposition by mathematical induction, for $n\ge2$.  
\par For $n=2$, $S(n)$ is true since we know from \cite{1}, \cite{2}, that the standard $2$-inner product, $(x,y|z)_*=\det\begin{pmatrix} \langle x,y \rangle & \langle x,z \rangle \\ \langle z,y \rangle &\langle z,z \rangle \end{pmatrix}$, satisfies conditions $P1)-P6)$.
\par Suppose $S(n)$ is true and prove that proposition $S(n+1)$ is true. The $(n+1)$-iterated $2$-inner product is given by
$$(x,y|x_{n+1},x_{n},...,x_{2})_*=\begin{vmatrix}
(x,y|x_{n},...,x_{2})_* & (x,x_{n+1}|x_{n},...,x_{2})_*\\ (x_{n+1},y|x_{n},...,x_{2})_* &(x_{n+1},x_{n+1}|x_{n},...,x_{2})_*
\end{vmatrix}.$$

\par Let us prove P1) for $n+1$. First we prove that $(x,x|x_{n+1},...,x_{2})_* \geq 0$, for $x, x_2,\ldots x_{n+1}\in X$.
 \par\underline{Case 1:} $(x,x|x_{n},...,x_{2})_*=0 $. Then, from property P1) for $n$, it results that $x, x_{2},...,x_{n}$ are linearly dependent. 
From the hypothesis of induction and from Lemma \ref{L1a} it follows that $(x,x_{n+1}|x_{n},...,x_{2})_*=0$. Then
\begin{eqnarray*}
(x,x|x_{n+1},\ldots,x_{2})_*&=& \begin{vmatrix} (x,x|x_{n},...,x_{2})_*
& (x,x_{n+1}|x_{n},...,x_{2})_*\\ (x_{n+1},x|x_{n},...,x_{2})_* &(x_{n+1},x_{n+1}|x_{n},...,x_{2})_*
\end{vmatrix}\\
&=&\begin{vmatrix}
0 &0\\ (x_{n+1},x|x_{n},...,x_{2})_* &(x_{n+1},x_{n+1}|x_{n},...,x_{2})_*
\end{vmatrix}=0.
\end{eqnarray*}
 \par\underline{Case 2:} $(x,x|x_{n},...,x_{2})_*>0 $. From P1) for $n$ we have  
$(z,z|x_{n},...,x_{2})_* \geq 0$, for all $z\in X$, then 
$$(\lambda x+x_{n+1},\lambda x +x_{n+1}|x_{n},...,x_{2})_* \geq 0,\;\mbox{ for all}\;\lambda \in \mathbb{R}.$$
We obtain the following relation:
$$\lambda^{2}(x,x|x_{n},...,x_{2})_*+2\lambda(x,x_{n+1}|x_{n},...,x_{2})_*+(x_{n+1},x_{n+1}|x_{n},...,x_{2})_* \geq 0,\;\forall \lambda.$$
Since $(x,x|x_{n},...,x_{2})_*>0$, the discriminant $\Delta_{\lambda}$ of this polynomial in variable $\lambda$ is not strictly positive. Hence
$(x,x|x_{n+1},x_{n},...,x_{2})_* = -\frac14\Delta_{\lambda}\ge 0.$ 
So, in both cases we obtain $(x,x|x_{n+1},...,x_{2})_* \geq 0$. 

\par On the other hand, let us suppose that $(x,x|x_{n+1},x_n,...,x_2)_*=0$, which means that 
$$(x,x|x_n,...,x_2)_*(x_{n+1},x_{n+1}|x_n,...,x_2)_*-(x,x_{n+1}|x_n,...,x_2)_*^2=0.$$  
If $(x_{n+1},x_{n+1}|x_n,...,x_2)_*\not=0$, the expression above is equal to  $-\frac14\Delta_{\lambda}$, where $\Delta_{\lambda}$ is the discriminat of the polynomial equation of degree 2 in $\lambda$: $Q(\lambda)=0$, where $Q(\lambda)=(x+\lambda x_{n+1},x+\lambda x_{n+1}|x_n,...,x_2)_*$. Since the discriminat is $0$, then there exists $\lambda_0\in\mathbb{R}$, for which $Q(\lambda_0)=0$. From condition P1) for $n$ it follows that $x+\lambda_0x_{n+1}, x_n,...,x_2$ are linearly dependent. Then, there are the numbers $\alpha, \alpha_i\in\mathbb{R}$, not all null, such that $\alpha(x+\lambda_0x_{n+1})+\sum_{i=2}^{n}\alpha_i x_i=0$. Therefore, $x, x_2,...,x_{n+1}$ are linearly dependent. If $(x_{n+1},x_{n+1}|x_n,...,x_2)_*=0$, then $x_2,...,x_{n+1}$ are linearly dependent from P1) for $n$. Then $x, x_2,...,x_{n+1}$ are linearly dependent. Condition P1) is completely proved for $n+1$.

\par  We prove condition P2) for $n+1$: 
$$(x,x|x_{n+1},...,x_{2})_*=\begin{vmatrix}
(x,x|x_{n},...,x_{2})_* & (x,x_{n+1}|x_{n},...,x_{2})_*\\ (x_{n+1},x|x_{n},...,x_{2})_* &(x_{n+1},x_{n+1}|x_{n},...,x_{2})_*
\end{vmatrix} $$
$$=\begin{vmatrix}
(x_{n+1},x_{n+1}|x_{n},...,x_{2})_* & (x_{n+1},x|x_{n},...,x_{2})_*\\ (x,x_{n+1}|x_{n},...,x_{2})_* &(x,x|x_{n},...,x_{2})_*
\end{vmatrix} =(x_{n+1},x_{n+1}|x,x_{n},...,x_{2})_*.$$
Consequently, condition P2) is true for $n+1$. 
\par We pass to the verification of condition P3). We have
$$(x,y|x_{n+1},x_{n},...,x_{2})_*=\begin{vmatrix}
(x,y|x_{n},...,x_{2})_* & (x,x_{n+1}|x_{n},...,x_{2})_*\\ (x_{n+1},y|x_{n},...,x_{2})_* &(x_{n+1},x_{n+1}|x_{n},...,x_{2})_*
\end{vmatrix} $$
$$ = \begin{vmatrix}
(y,x|x_{n},...,x_{2})_* & (x_{n+1},y|x_{n},...,x_{2})_*\\ (x,x_{n+1}|x_{n},...,x_{2})_* &(x_{n+1},x_{n+1}|x_{n},...,x_{2})_*
\end{vmatrix}=(y,x|x_{n+1},x_{n},...,x_{2})_*, $$
because $\det A=\det A^{T}$, for any square matrix $A$ and $(x,y|x_{n},...,x_{2})_*=(y,x|x_{n},..,x_{2})_*$. So, the $(n+1)$-iterated $2$-inner product satisfies condition  P3) for $n+1$. 
\par We pass now to condition P4). Since we have
$$(\alpha x,y|x_{n+1},x_{n},...,x_{2})_*= \begin{vmatrix}
(\alpha x,y|x_{n},...,x_{2})_* & (\alpha x,x_{n+1}|x_{n},...,x_{2})_*\\ (x_{n+1},y|x_{n},...,x_{2})_* &(x_{n+1},x_{n+1}|x_{n},...,x_{2})_*
\end{vmatrix}$$
$$=\begin{vmatrix}
\alpha(x,y|x_{n},...,x_{2})_* & \alpha(x,x_{n+1}|x_{n},...,x_{2})_*\\ (x_{n+1},y|x_{n},...,x_{2})_* &(x_{n+1},x_{n+1}|x_{n},...,x_{2})_*
\end{vmatrix}=\alpha(x,y|x_{n+1},x_{n},...,x_{2})_*, $$
it follows that condition P4) is proved for $n+1$.
\par  For condition P5) for $n+1$, we take into account that $(x+x',y|x_{n+1},x_{n},...,x_{2})_*$ can be expressed by a determinat of second order, having on the first line the elements $(x+x',y|x_{n},...,x_{2})_*$ and $(x+x',x_{n+1}|x_{n},...,x_{2})_*$, respectively, and on the second line, elements which do not depend on $x$ and $x'$. Then, using by induction the additivity in the first  argument of the products above, and then the additivity of the determinant with regard to the first line, it follows immediately that  
\begin{eqnarray*}
(x+x',y|x_{n+1},x_{n},...,x_{2})_*=(x,y|x_{n+1},x_{n},...,x_{2})_*+(x',y|x_{n+1},x_{n},...,x_{2})_*.
\end{eqnarray*}

 \framebox

\begin{proposition}
Let $(X,\langle \cdot,\cdot\rangle)$ be a real pre-Hilbert space. For $x,y,x_2,\ldots,x_n\in X$, $n\ge 2$ and $t\in\mathbb{R}$:
\begin{equation}
(tx,ty|tx_2,\ldots,tx_n)_*=t^{2^{n}}(x,y|x_2,\ldots, x_n)_*.
\end{equation}
\end{proposition}
\par{\it Proof} For $n=2$, $(tx,ty|tx_2)_*=t^4(x,y|x_2)_*.$ Then, it follows by  mathematical induction. \framebox

\begin{remark}\label{R2.4}{\rm Theorem \ref{T3} allows us to furnish an example of weak $n$-inner product which is not a $n$-inner product. For this, let $X=\mathbb{R}^3$ endowed with the usual inner product. Then, from Theorem \ref{T3}, $3$-iterated $2$-inner product $(\cdot,\cdot|\cdot,\cdot)_*:X^4\to\mathbb{R}$ is a weak $3$-inner product, but it is not a $3$-inner product. Indeed, if axiom I2) would be true for $3$-iterated $2$-inner product then we must have:
\begin{equation}\label{e2.7}
(x,x|u,v)_*=(v,v|u,x)_*,\quad\mbox{for all}\;x,u,v\in X.
\end{equation} 
But, if we choose $x=(1,0,0)$, $u=(1,1,1)$, $v=(2,1,2)$ we have
\begin{eqnarray*}
(x,x|u,v)_*&=&
\begin{vmatrix}\;\;\begin{vmatrix}\langle x,x\rangle& \langle x,v\rangle\\
\langle v,x \rangle&\langle v,v\rangle\end{vmatrix}\;
\begin{vmatrix}\langle x,u\rangle& \langle x,v\rangle\\
\langle v,u \rangle&\langle v,v\rangle\end{vmatrix}\\
\;&\;\\
\;\;\begin{vmatrix}\langle u,x\rangle& \langle u,v\rangle\\
\langle v,x \rangle&\langle v,v\rangle\end{vmatrix}\;
\begin{vmatrix}\langle u,u\rangle& \langle u,v\rangle\\
\langle v,u \rangle&\langle v,v\rangle\end{vmatrix}
\end{vmatrix}=
\begin{vmatrix}\;\begin{vmatrix}1& 2\\
2&9\end{vmatrix}\;
\begin{vmatrix}1&2\\
5&9\end{vmatrix}\\
\;&\;\\
\;\begin{vmatrix}1& 5\\
2&9\end{vmatrix}\;
\begin{vmatrix}3& 5\\
5&9\end{vmatrix}
\end{vmatrix}
=\begin{vmatrix}5&-1\\-1&2\end{vmatrix}=9
\end{eqnarray*} 
and on other hand
\begin{eqnarray*}
(v,v|u,x)_*&=&
\begin{vmatrix}\;\;\begin{vmatrix}\langle v,v\rangle& \langle v,x\rangle\\
\langle x,v \rangle&\langle x,x\rangle\end{vmatrix}\;
\begin{vmatrix}\langle v,u\rangle& \langle v,x\rangle\\
\langle x,u \rangle&\langle x,x\rangle\end{vmatrix}\\
\;&\;\\
\;\;\begin{vmatrix}\langle u,v\rangle& \langle u,x\rangle\\
\langle x,v \rangle&\langle x,x\rangle\end{vmatrix}\;
\begin{vmatrix}\langle u,u\rangle& \langle u,x\rangle\\
\langle x,u \rangle&\langle x,x\rangle\end{vmatrix}
\end{vmatrix}=
\begin{vmatrix}\;\begin{vmatrix}9& 2\\
2&1\end{vmatrix}\;
\begin{vmatrix}5&2\\
1&1\end{vmatrix}\\
\;&\;\\
\;\begin{vmatrix}5&1\\
2&1\end{vmatrix}\;
\begin{vmatrix}3& 1\\
1&1\end{vmatrix}
\end{vmatrix}
=\begin{vmatrix}5&3\\3&2\end{vmatrix}=1.
\end{eqnarray*} 
Hence relation (\ref{e2.7}) is not true. Consequently axiom I2) is not satisfied. Therefore $3$-iterated $2$-inner product is not a $3$-inner product.
}\end{remark}

\section{Representation of the  $n$-iterated $2$-inner product in terms of the standard $k$-inner products, $(k\le n)$}

\par In this section we obtain a representation of the $n$-iterated $2$-inner product, given in Definition \ref{D3} in terms of the standard $k$-inner products $k\le n$. For this we use Dodgson's identity for determinants, \cite{dodgson1}, \cite{dodgson2}. Historical notes about this identity, in connection with Chi\`o's formula can be found in \cite{abeles}. To express this identity we adopt the compact notation used by Eves \cite{eves}. If $A=(a_{i,j})_{1\le i,j\le n}$ is a square matrix, denote the determinant of $A$ by $|a_{1,1}\;\ldots a_{n,n}|$ and the sub-determinant involving rows $i_1,\ldots,i_s$ and columns $j_1,\ldots, j_s$ by $|a_{i_1,j_1}\;\ldots a_{i_s,j_s}|$. In \cite{eves} - Theorem 3.6.3, the following Dodgson type identity $(n\ge3)$ is proved:
\begin{eqnarray}
&&|a_{1,1}\;\ldots\;a_{n-2,n-2}|\cdot|a_{1,1}\;\ldots\;a_{n,n}|\nonumber\\
&&\quad\qquad=\begin{vmatrix} 
|a_{1,1}\;\ldots\;a_{n-2,n-2}\; a_{n-1,n-1}|&|a_{1,1}\;\ldots a_{n-2.n-2}\;a_{n-1,n}|\\
|a_{1,1}\;\ldots\;a_{n-2,n-2}\; a_{n,n-1}|&|a_{1,1}\;\ldots a_{n-2,n-2}\;a_{n,n}|
\end{vmatrix}.\label{e40}
\end{eqnarray}

For us it is more convenient to use the following identity $(n\ge 3)$:
\begin{eqnarray}
&&|a_{2,2}\;\ldots\;a_{n-1,n-1}|\cdot|a_{1,1}\;\ldots\;a_{n,n}|\nonumber\\
&&\quad\qquad=\begin{vmatrix} 
\;|a_{1,1}\;\ldots\;a_{n-2,n-2}\; a_{n-1,n-1}|&|a_{1,2}\;\ldots a_{n-2,n-1}\;a_{n-1,n}|\;\\
\;|a_{2,1}\;\ldots\;a_{n-1,n-2}\; a_{n,n-1}|&|a_{2,2}\;\ldots a_{n-1,n-1}\;a_{n,n}|\;
\end{vmatrix}.\label{e41}
\end{eqnarray}

For $n=3$ one has:

\begin{equation}\label{e41a}
a_{2,2}\cdot |a_{1,1}\;a_{2,2}\;a_{3,3}|=\begin{vmatrix} \;|a_{1,1}\;a_{2,2}|&|a_{1,2}\;a_{2,3}|\;\\
\;|a_{2,1}\;a_{3,2}|&|a_{2,2}\;a_{3,3}|\; \end{vmatrix}.
\end{equation}   

Formula (\ref{e41}) can be easily obtained applying formula (\ref{e40}). Indeed, first note that
$$|a_{1,1}\;\ldots\;a_{n,n}|=(-1)^{n-1}|a_{1,2}\;\ldots\;a_{n-1,n}\;a_{n,1}|=|a_{2,2}\;\ldots\;a_{n,n}\;a_{1,1}|.$$ 
Then, applying rule (\ref{e40}) for our new matrix we find
\begin{eqnarray*}
&&|a_{2,2}\;\ldots \;a_{n-1,n-1}|\cdot|a_{2,2}\;\ldots\;a_{n,n}\;a_{1,1}|\\
&=&\begin{vmatrix} 
|a_{2,2}\;\ldots\;a_{n-1,n-1}\;a_{n,n}|&|a_{2,2}\;\ldots a_{n-1,n-1}\;a_{n,1}|\;\\
|a_{2,2}\;\ldots\;a_{n-1,n-1}\;a_{1,n}|&|a_{2,2}\;\ldots a_{n-1,n-1}\;a_{1,1}|\;
\end{vmatrix}\\
&=&\begin{vmatrix} 
|a_{2,2}\;\ldots\;a_{n-1,n-1}\;a_{n,n}|&(-1)^{n-2}|a_{2,1}\;\ldots a_{n-1,n-2}\;a_{n,n-1}|\;\\
(-1)^{n-2}|a_{1,2}\;\ldots\;a_{n-2,n-1}\;a_{n-1,n}|&(-1)^{2n-4}|a_{1,1}\;\ldots a_{n-2,n-2}\;a_{n-1,n-1}|\;
\end{vmatrix}\\
&=&\begin{vmatrix} 
\;|a_{1,1}\;\ldots\;a_{n-2,n-2}\; a_{n-1,n-1}|&|a_{1,2}\;\ldots a_{n-2,n-1}\;a_{n-1,n}|\;\\
\;|a_{2,1}\;\ldots\;a_{n-1,n-2}\; a_{n,n-1}|&|a_{2,2}\;\ldots a_{n-1,n-1}\;a_{n,n}|\;
\end{vmatrix}.
\end{eqnarray*}

\par Note that, conversely, from relation (\ref{e41}) one can deduce relation (\ref{e40}).

\par Let $(X,\langle\cdot,\cdot\rangle)$ be an inner product space. For $x,y,z,v\in X$, from (\ref{e41}), for $n=3$ we deduce

\begin{eqnarray*}
(x,y|v,z)_*&=&\begin{vmatrix} (x,y|z)_*&(x,v|z)_*\\(v,y|z)_*&(v,v|z)_*\end{vmatrix} 
=\begin{vmatrix}\;\begin{vmatrix} \langle x,y\rangle&\langle x,z\rangle\\ \langle z,y\rangle&\langle z,z\rangle\end{vmatrix} \; 
\begin{vmatrix} \langle x,v\rangle&\langle x,z\rangle\\\langle z,v\rangle&\langle z,z\rangle\end{vmatrix}\;\\
\quad\\
\;\begin{vmatrix} \langle v,y\rangle&\langle v,z\rangle\\\langle z,y\rangle&\langle z,z\rangle\end{vmatrix}\;  \begin{vmatrix} \langle v,v\rangle&\langle v,z\rangle\\\langle z,v\rangle&\langle z,z\rangle\end{vmatrix}\;\end{vmatrix}\\
&=&\begin{vmatrix}\;\begin{vmatrix} \langle x,y\rangle&\langle x,z\rangle\\ \langle z,y\rangle&\langle z,z\rangle\end{vmatrix} \; 
\begin{vmatrix}   \langle x,z\rangle&\langle x,v\rangle\\ \langle z,z\rangle&\langle z,v\rangle\end{vmatrix}\;\\
\quad\\
\;\begin{vmatrix} \langle z,y\rangle&\langle z,z\rangle\\ \langle v,y\rangle&\langle v,z\rangle\end{vmatrix}\;  
  \begin{vmatrix} \langle z,z\rangle&\langle z,v\rangle\\ \langle v,z\rangle&\langle v,v\rangle\end{vmatrix}\;\end{vmatrix}
=\langle z,z\rangle \begin{vmatrix}\langle x,y\rangle&\langle x,z\rangle&\langle x,v\rangle\\
\langle z,y\rangle&\langle z,z\rangle&\langle z,v\rangle\\
\langle v,y\rangle&\langle v,z\rangle&\langle v,v\rangle
\end{vmatrix}\\
&=&\langle z,z\rangle 
 \begin{vmatrix}\langle x,y\rangle&\langle x,v\rangle&\langle x,z\rangle\\
\langle v,y\rangle&\langle v,v\rangle&\langle v,z\rangle\\
\langle z,y\rangle&\langle z,v\rangle&\langle z,z\rangle\end{vmatrix}.
\end{eqnarray*}
Hence, we obtained
\begin{equation}\label{e45a}
(x,y|v,z)_*=\langle z,z\rangle \langle x,y|v,z\rangle.
\end{equation}

\par Also, using formula (\ref{e41}), for $n=4$ and then formula (\ref{e45a}) we obtain : 

\begin{eqnarray*}
&&\langle z,z|w\rangle\langle z,z\rangle^2\begin{vmatrix} \langle x,y\rangle&\langle x,z\rangle&\langle x,v\rangle&\langle x,w\rangle\\
\langle z,y\rangle&\langle z,z\rangle&\langle z,v\rangle&\langle z,w\rangle\\
\langle v,y\rangle&\langle v,z\rangle&\langle v,v\rangle&\langle v,w\rangle\\
\langle w,y\rangle&\langle w,z\rangle&\langle w,v\rangle&\langle w,w\rangle\end{vmatrix}\\
&=&\langle z,z\rangle^2 \begin{vmatrix}\langle z,z\rangle&\langle z,w\rangle\\ \langle w,z\rangle&\langle w,w\rangle\end{vmatrix}
\begin{vmatrix} \langle x,y\rangle&\langle x,z\rangle&\langle x,w\rangle&\langle x,v\rangle\\
\langle z,y\rangle&\langle z,z\rangle&\langle z,w\rangle&\langle z,v\rangle\\
\langle w,y\rangle&\langle w,z\rangle&\langle w,w\rangle&\langle w,v\rangle\\
\langle v,y\rangle&\langle v,z\rangle&\langle v,w\rangle&\langle v,v\rangle\end{vmatrix}\\
\end{eqnarray*}
\begin{eqnarray*}
&=&\langle z,z\rangle^2 
\begin{vmatrix}\;\begin{vmatrix} \langle x,y\rangle&\langle x,z\rangle&\langle x,w\rangle\\
\langle z,y\rangle&\langle z,z\rangle&\langle z,w\rangle\\
\langle w,y\rangle&\langle w,z\rangle&\langle w,w\rangle\end{vmatrix}\;
\begin{vmatrix}\langle x,z\rangle&\langle x,w\rangle&\langle x,v\rangle\\
\langle z,z\rangle&\langle z,w\rangle&\langle z,v\rangle\\
\langle w,z\rangle&\langle w,w\rangle&\langle w,v\rangle\end{vmatrix}\\
\quad\\
\;\begin{vmatrix}\langle z,y\rangle&\langle z,z\rangle&\langle z,w\rangle\\
\langle w,y\rangle&\langle w,z\rangle&\langle w,w\rangle\\
\langle v,y\rangle&\langle v,z\rangle&\langle v,w\rangle\end{vmatrix}\;
\begin{vmatrix}\langle z,z\rangle&\langle z,w\rangle&\langle z,v\rangle\\
\langle w,z\rangle&\langle w,w\rangle&\langle w,v\rangle\\
\langle v,z\rangle&\langle v,w\rangle&\langle v,v\rangle\end{vmatrix}\;
\end{vmatrix}\\ 
&=&
\begin{vmatrix}\;\langle z,z\rangle\begin{vmatrix} \langle x,y\rangle&\langle x,w\rangle&\langle x,z\rangle\\
\langle w,y\rangle&\langle w,w\rangle&\langle w,z\rangle\\
\langle z,y\rangle&\langle z,w\rangle&\langle z,z\rangle\end{vmatrix}\;
\langle z,z\rangle\begin{vmatrix}\langle x,v\rangle&\langle x,w\rangle&\langle x,z\rangle\\
\langle w,v\rangle&\langle w,w\rangle&\langle w,z\rangle\\
\langle z,v\rangle&\langle z,w\rangle&\langle z,z\rangle\end{vmatrix}\;\\
\quad\\
\;\langle z,z\rangle\begin{vmatrix}\langle v,y\rangle&\langle v,w\rangle&\langle v,z\rangle\\
\langle w,y\rangle&\langle w,w\rangle&\langle w,z\rangle\\
\langle z,y\rangle&\langle z,w\rangle&\langle z,z\rangle\end{vmatrix}\;
\;\langle z,z\rangle\begin{vmatrix}\langle v,v\rangle&\langle v,w\rangle&\langle v,z\rangle\\
\langle w,v\rangle&\langle w,w\rangle&\langle w,z\rangle\\
\langle z,v\rangle&\langle z,w\rangle&\langle z,z\rangle\end{vmatrix}\;
\end{vmatrix}\\
&=&\begin{vmatrix} \langle z,z\rangle\langle x,y|w,z\rangle& \langle z,z\rangle \langle x,v|w,z\rangle\\ \langle z,z\rangle \langle v,y|w,z\rangle& \langle z,z\rangle\langle v,v|w,z\rangle\end{vmatrix}\\
&=&\begin{vmatrix}  (x,y|w,z)_*&(x,v|w,z)_*\\
(v,y|w,z)_*&(v,v|w,z)_*\end{vmatrix}\\
&=&(x,y|v,w,z)_*.
\end{eqnarray*}

Hence
\begin{equation}\label{e45.1}
(x,y|v,w,z)_*=\langle z,z|w\rangle\langle z,z\rangle^2\langle x,y|v,w,z\rangle.
\end{equation}

\par The results obtained in (\ref{e45a}) and (\ref{e45.1}) can be generalized as it is shown in the next theorem. We extend the definition of the standard weak $n$-inner product, for $n=1$,  by the convention $ \langle x,y|x_1,\ldots,x_2\rangle=\langle x,y\rangle$.

\begin{theorem}\label{T4}
Let $(X,\langle\cdot,\cdot\rangle)$ be an inner product space. For $n\ge 2$, consider the vectors $x,y,x_2,\ldots,x_n\in X.$. 
Then 
\begin{equation}\label{e46}
(x,y|x_n,\ldots,x_2)_*=E_n\cdot\langle x,y|x_n,\ldots,x_2\rangle,
\end{equation}
where $E_2=1$ and
\begin{equation}
E_n=\prod_{k=2}^{n-1}\langle x_k,x_k|x_{k-1},\ldots,x_2\rangle^{2^{n-k-1}},\; (n\ge 3).
\end{equation}
\end{theorem}
\par{\it Proof}\ For $n=2$ the theorem is immediate, since $(x,y|x_2)_*=\langle x,y|x_2\rangle$ and $E_2=1$. For $n=3$  the theorem follows from relation (\ref{e45a}), for the choice $v=x_3$ and $z=x_2$. Then $E_3=\langle z,z\rangle$. For $n\ge 4$ we prove by induction. Suppose the theorem true for $n\ge 3$ and let us prove it for $n+1$.
Using the hypothesis of induction we get 
\begin{eqnarray}
&&(x,y|x_{n+1},x_n,\ldots,x_2)_*=\begin{vmatrix} (x,y|x_n,\ldots, x_2)_*&(x,x_{n+1}|x_n,\ldots, x_2)_*\\(x_{n+1},y|x_n,\ldots, x_2)_*&(x_{n+1},x_{n+1}|x_n,\ldots, x_2)_*\end{vmatrix}\nonumber\\ 
&&\qquad\qquad=\begin{vmatrix} E_n\langle x,y|x_n,\ldots, x_2\rangle &E_n\langle x,x_{n+1}|x_n,\ldots, x_2\rangle\\E_n\langle x_{n+1},y|x_n,\ldots, x_2\rangle &E_n\langle x_{n+1},x_{n+1}|x_n,\ldots, x_2\rangle \end{vmatrix}\nonumber\\ 
&&\qquad\qquad=(E_n)^2\begin{vmatrix} \langle x,y|x_n,\ldots, x_2\rangle &\langle x,x_{n+1}|x_n,\ldots, x_2\rangle\\\langle x_{n+1},y|x_n,\ldots, x_2\rangle &\langle x_{n+1},x_{n+1}|x_n,\ldots, x_2\rangle\end{vmatrix}.
\label{e48}
\end{eqnarray}
We transform all the four elements from the above determinant. Each of them is a determinant of order $n$. First, in the following determinant, changing the order of the last $n-1$ lines and then changing the order of the last $n-1$ columns we obtain successively
\begin{eqnarray} 
\langle x,y|x_n,\ldots, x_2\rangle&=&
\begin{vmatrix}
\langle x,y\rangle&\langle x,x_n\rangle&\dots&\langle x,x_2\rangle\\
\langle x_n,y\rangle&\langle x_n,x_n\rangle&\dots&\langle x_n,x_2\rangle\\
\vdots&\vdots&\ddots&\vdots\\
\langle x_2,y\rangle&\langle x_2,x_n\rangle&\dots&\langle x_2,x_2\rangle
\end{vmatrix}\nonumber\\
&=&(-1)^{\frac{(n-1)(n-2)}2}
\begin{vmatrix}
\langle x,y\rangle&\langle x,x_n\rangle&\dots&\langle x,x_2\rangle\\
\langle x_2,y\rangle&\langle x_2,x_n\rangle&\dots&\langle x_2,x_2\rangle\\
\vdots&\vdots&\ddots&\vdots\\
\langle x_n,y\rangle&\langle x_n,x_n\rangle&\dots&\langle x_n,x_2\rangle
\end{vmatrix}\nonumber\\
&=&\begin{vmatrix}
\langle x,y\rangle&\langle x,x_2\rangle&\dots&\langle x,x_n\rangle\\
\langle x_2,y\rangle&\langle x_2,x_2\rangle&\dots&\langle x_2,x_n\rangle\\
\vdots&\vdots&\ddots&\vdots\\
\langle x_n,y\rangle&\langle x_n,x_2\rangle&\dots&\langle x_n,x_n\rangle
\end{vmatrix}.\label{e50}
\end{eqnarray}
Next, for the second determinant, we change the order of all the $n$ columns and then we change the order of the last $n-1$ lines we obtain:

\begin{eqnarray} 
 \langle x,x_{n+1}|x_n,\ldots, x_2\rangle&=&
\begin{vmatrix}
\langle x,x_{n+1}\rangle&\langle x,x_n\rangle&\dots&\langle x,x_2\rangle\\
\langle x_n,x_{n+1}\rangle&\langle x_n,x_n\rangle&\dots&\langle x_n,x_2\rangle\\
\vdots&\vdots&\ddots&\vdots\\
\langle x_2,x_{n+1}\rangle&\langle x_2,x_n\rangle&\dots&\langle x_2,x_2\rangle
\end{vmatrix}\nonumber\\
&=&(-1)^{\frac{(n-1)n}2}
\begin{vmatrix}
\langle x,x_2\rangle&\langle x,x_3\rangle&\dots&\langle x,x_{n+1}\rangle\\
\langle x_n,x_2\rangle&\langle x_n,x_3\rangle&\dots&\langle x_n,x_{n+1}\rangle\\
\vdots&\vdots&\ddots&\vdots\\
\langle x_2,x_2\rangle&\langle x_2,x_3\rangle&\dots&\langle x_2,x_{n+1}\rangle
\end{vmatrix}\nonumber\\
\quad\nonumber\\
&=&(-1)^{(n-1)^2}\begin{vmatrix}
\langle x,x_2\rangle&\langle x,x_3\rangle&\dots&\langle x,x_{n+1}\rangle\\
\langle x_2,x_2\rangle&\langle x_2,x_3\rangle&\dots&\langle x_2,x_{n+1}\rangle\\
\vdots&\vdots&\ddots&\vdots\\
\langle x_n,x_{n+1}\rangle&\langle x_n,x_3\rangle&\dots&\langle x_n,x_{n+1}\rangle\\
\end{vmatrix},\label{e51}
\end{eqnarray}
since $\frac{(n-1)n}2+\frac{(n-2)(n-1)}2=(n-1)^2.$

Similar operations there can be made for the third determinant. We change the order of all the $n$ lines and then we change the order of the last $n-1$ columns and we get:

\begin{eqnarray} 
 \langle x_{n+1},y|x_n,\ldots, x_2\rangle&=&
\begin{vmatrix}
\langle x_{n+1},y\rangle&\langle x_{n+1},x_n\rangle&\dots&\langle x_{n+1},x_2\rangle\\
\langle x_n,y\rangle&\langle x_n,x_n\rangle&\dots&\langle x_n,x_2\rangle\\
\vdots&\vdots&\ddots&\vdots\\
\langle x_2,y\rangle&\langle x_2,x_n\rangle&\dots&\langle x_2,x_2\rangle
\end{vmatrix}\nonumber\\
&=&(-1)^{\frac{(n-1)n}2}
\begin{vmatrix}
\langle x_2,y\rangle&\langle x_2,x_n\rangle&\dots&\langle x_2,x_2\rangle\\
\langle x_3,y\rangle&\langle x_3,x_n\rangle&\dots&\langle x_3,x_2\rangle\\
\vdots&\vdots&\ddots&\vdots\\
\langle x_{n+1},y\rangle&\langle x_{n+1},x_n\rangle&\dots&\langle x_{n+1},x_2\rangle
\end{vmatrix}\nonumber\\
\quad\nonumber\\
&=&(-1)^{(n-1)^2}\begin{vmatrix}
\langle x_2,y\rangle&\langle x_2,x_2\rangle&\dots&\langle x_2,x_n\rangle\\
\langle x_3,y\rangle&\langle x_3,x_2\rangle&\dots&\langle x_3,x_n\rangle\\
\vdots&\vdots&\ddots&\vdots\\
\langle x_{n+1},y\rangle&\langle x_{n+1},x_2\rangle&\dots&\langle x_{n+1},x_n\rangle\\
\end{vmatrix}.\label{e52}
\end{eqnarray}

Finally, applying formula I2) we have
\begin{eqnarray} 
 \langle x_{n+1},x_{n+1}|x_n,\ldots, x_2\rangle&=& \langle x_2,x_2|x_3,\ldots,x_{n+1}\rangle\nonumber\\
&=&
\begin{vmatrix}
\langle x_2,x_2\rangle&\langle x_2,x_3\rangle&\dots&\langle x_2, x_{n+1}\rangle\\
\langle x_3,x_2\rangle&\langle x_3,x_3\rangle&\dots&\langle  x_3, x_{n+1}\rangle\\
\vdots&\vdots&\ddots&\vdots\\
\langle x_{n+1},x_2\rangle&\langle x_{n+1},x_3\rangle&\dots&\langle x_{n+1},x_{n+1}\rangle
\end{vmatrix}.\label{e53}
\end{eqnarray}

Consider the matrix
$$A=\begin{vmatrix}
\langle x,y\rangle&\langle x,x_2\rangle&\dots&\langle x, x_{n+1}\rangle\\
\langle x_2,y\rangle&\langle x_2,x_2\rangle&\dots&\langle  x_2, x_{n+1}\rangle\\
\vdots&\vdots&\ddots&\vdots\\
\langle x_{n+1},y\rangle&\langle x_{n+1},x_2\rangle&\dots&\langle x_{n+1},x_{n+1}\rangle
\end{vmatrix}.$$

Denote the elements of $A$ by $a_{i,j}$, $1\le i,j\le n+1$. Using the notation given in the beginning of the section we can write $|A|=|a_{1,1}\;a_{2,2}\;\ldots\;a_{n+1,n+1}|$.  

\par From formula (\ref{e50}) we obtain $\langle x,y|x_n,\ldots,x_2\rangle=|a_{1,1}\;a_{2,2}\;\ldots\;a_{n,n}|$.
\par From formula (\ref{e51}) we obtain $\langle x,x_{n+1}|x_n,\ldots,x_2\rangle=(-1)^{(n-1)^2}|a_{1,2}\;a_{2,3}\;\ldots\;a_{n,n+1}|$.
\par From formula (\ref{e52}) we obtain $\langle x_{n+1},y|x_n,\ldots,x_2\rangle=(-1)^{(n-1)^2}|a_{2,1}\;a_{3,2}\;\ldots\;a_{n+1,n}|$.
\par From formula (\ref{e53}) we obtain $\langle x_{n+1},x_{n+1}|x_n,\ldots,x_2\rangle=|a_{2,2}\;a_{3,3}\;\ldots\;a_{n+1,n+1}|$.

Then applying formula (\ref{e41}) for $n+1$ instead of $n$ we arrive to

\begin{eqnarray}
&&\begin{vmatrix} \langle x,y|x_n,\ldots, x_2\rangle &\langle x,x_{n+1}|x_n,\ldots, x_2\rangle\\\langle x_{n+1},y|x_n,\ldots, x_2\rangle&\langle x_{n+1},x_{n+1}|x_n,\ldots, x_2\rangle\end{vmatrix}\nonumber\\
&=&
|a_{1,1}\;a_{2,2}\;\ldots\;a_{n+1,n+1}|\cdot |a_{2,2}\;a_{3,3}\;\ldots\;a_{n,n}|.\label{e60}
\end{eqnarray}

If we change the order of the last $n$ lines and  of the last $n$ columns  in  $|A|$ the determinant does not change, i.e.
$$|a_{1,1}\;a_{2,2}\;\ldots\;a_{n+1,n+1}|=|a_{1,1}\;a_{n+1,n+1}\;a_{n,n}\;\ldots\;a_{2,2}|.$$

But 
\begin{eqnarray*}
|a_{1,1}\;a_{n+1,n+1}\;a_{n,n}\;\ldots\;a_{2,2}|
&=&\begin{vmatrix}
\langle x,y\rangle&\langle x,x_{n+1}\rangle&\dots&\langle x, x_2\rangle\\
\langle x_{n+1},y\rangle&\langle x_{n+1},x_{n+1}\rangle&\dots&\langle  x_{n+1},x_2\rangle\\
\vdots&\vdots&\ddots&\vdots\\
\langle x_2,y\rangle&\langle x_2,x_{n+1}\rangle&\dots&\langle x_2,x_2\rangle
\end{vmatrix}\\
&=&\langle x,y|x_{n+1},x_n,\ldots, x_2\rangle.
\end{eqnarray*}
Therefore,
\begin{equation}\label{e62}
|a_{1,1}\;a_{2,2}\;\ldots\;a_{n+1,n+1}|=\langle x,y|x_{n+1},x_n,\ldots, x_2\rangle.
\end{equation}

Also, if we change the order of the lines and columns in determinat  $|a_{2,2}\;a_{3,3}\;\ldots\;a_{n,n}|$, the value does not change. Hence
$$|a_{2,2}\;a_{3,3}\;\ldots\;a_{n,n}|=|a_{n,n}\;a_{n-1,n-1}\; \ldots\; a_{2,2}|.$$

But
\begin{eqnarray*}
|a_{n,n}\;a_{n-1,n-1}\; \ldots\; a_{2,2}|
&=&\begin{vmatrix}
\langle x_n,x_n\rangle&\langle x_n,x_{n-1}\rangle&\dots&\langle x_n, x_2\rangle\\
\langle x_{n-1},x_n\rangle&\langle x_{n-1},x_{n-1}\rangle&\dots&\langle  x_{n-1},x_2\rangle\\
\vdots&\vdots&\ddots&\vdots\\
\langle x_2,x_n\rangle&\langle x_2,x_{n-1}\rangle&\dots&\langle x_2,x_2\rangle
\end{vmatrix}\\
&=&(x_n,x_n|x_{n-1},\ldots,x_2).
\end{eqnarray*}
Therefore,
\begin{equation}\label{e63}
|a_{2,2}\;a_{3,3}\;\ldots\;a_{n,n}|=\langle x_n,x_n|x_{n-1},\ldots,x_2\rangle.
\end{equation}

From relations (\ref{e60}), (\ref{e62}) and (\ref{e63}) we conclude that

\begin{eqnarray}
&&\begin{vmatrix} \langle x,y|x_n,\ldots, x_2\rangle&\langle x,x_{n+1}|x_n,\ldots, x_2\rangle\\\langle x_{n+1},y|x_n,\ldots, x_2\rangle&\langle x_{n+1},x_{n+1}|x_n,\ldots, x_2\rangle\end{vmatrix}\nonumber\\
&&\qquad\qquad=\langle x,y|x_{n+1},x_n,\ldots, x_2\rangle\langle x_n,x_n|x_{n-1},\ldots,x_2\rangle. \label{64}
\end{eqnarray}

\par Replacing in (\ref{e48}) we obtain

$$(x,y|x_{n+1},x_n,\ldots,x_2)_*=(E_n)^2\langle x,y|x_{n+1},x_n,\ldots, x_2\rangle\langle x_n,x_n|x_{n-1},\ldots,x_2\rangle.$$
Since $(E_n)^2\langle x_n,x_n|x_{n-1},\ldots,x_2\rangle=E_{n+1}$ it results, finally, that
$$(x,y|x_{n+1},x_n,\ldots,x_2)_*=E_{n+1}\langle x,y|x_{n+1},x_n,\ldots, x_2\rangle.$$
\framebox

\section{Several applications of the $n$  iterated $2$-inner product}\label{sec:4}

{\bf 1.} Let $X=(X,\langle \cdot,\cdot\rangle)$ be an inner product space. Let $x,w,z\in X$. From Definition \ref{D3} we deduce 
\begin{eqnarray}
(x,x|w,z)_*&=&(x,x|z)_*(w,w|z)_*-(x,w|z)_*(w,x|z)_*\nonumber\\ 
&=&\Vert x|z \Vert^{2}\Vert w|z \Vert^{2}-(x,w|z)_*^{2}.\label{ee1}
\end{eqnarray}

Relation (\ref{ee1}) can be written as
\begin{eqnarray}
(x,x|w,z)_*&=&(\Vert x\Vert^{2}\Vert w\Vert^{2}\Vert z\Vert^{2}+2\langle w,z \rangle\langle z,x \rangle\langle x,w\rangle-\Vert x\Vert^{2}\langle w,z\rangle^{2}\nonumber\\
&&-\Vert w\Vert^{2}\langle z,x\rangle^{2}-\Vert z\Vert^{2}\langle x,w\rangle^{2})\Vert z\Vert^2.\label{ee2}
\end{eqnarray}
Since $(x,x|w,z)_* \geq 0$, then we obtain the inequality from Lupu and Schwarz \cite{10} given by the following:

\begin{equation} 
\Vert x\Vert^{2}\langle w,z\rangle^{2}+  \Vert w\Vert^{2}\langle z,x\rangle^{2}+\Vert z\Vert^{2}\langle x,w\rangle^{2} \leq \Vert x\Vert^{2}\Vert w\Vert^{2}\Vert z\Vert^{2}+2\langle w,z \rangle\langle z,x \rangle\langle x,w\rangle. 
\end{equation}
\par\

\par {\bf 2.} Formula (\ref{e45a}) can be written in the form

\begin{equation}\label{ee3}
(x,y|w,z)_*= \langle x,y|w,z\rangle\Vert z\Vert^{2}=\begin{vmatrix}\langle x,y\rangle & \langle x,w \rangle & \langle x,z\rangle \\ \langle w,y\rangle & \langle w,w \rangle & \langle w,z\rangle \\ \langle z,y \rangle & \langle z,w \rangle & \langle z,z\rangle
\end{vmatrix}\cdot\Vert z\Vert^{2}.
\end{equation}
Therefore, for $\Vert z\Vert\neq 1$, we have $(x,y|w,z)_*$ $\neq\langle x,y|w,z\rangle.$ Also, since in the case $x=y$, the determinant in (\ref{ee3}) is the Gram's determinant  $\Gamma(x,w,z)$, from relation  (\ref{ee3}) we can deduce:
\begin{equation}\label{ee4}
(x,x|w,z)_*= \Gamma(x,w,z)\cdot\Vert z\Vert^{2}.
\end{equation}
Since, also $(x,x|z,w)_*=\Gamma(x,z,w)\cdot\Vert w\Vert^{2}$ and  $\Gamma(x,w,z)=\Gamma(x,z,w)$ it results
\begin{equation}\label{ee6}
(x,x|z,w)_*\Vert z\Vert^{2}=(x,x|w,z)_*\Vert w\Vert^{2}.
\end{equation}
\par\ 

\par{\bf 3.} From Theorem \ref{T4} for $n=4$ we find that the $4$ iterated $2$-inner product  can be given in the following way:
\begin{eqnarray}
(x,y|v,w,z)_*&=&\Vert w|z\Vert^2\Vert z\Vert^4\begin{vmatrix}
\langle x,y\rangle & \langle x,z \rangle & \langle x,w\rangle &\langle x,v\rangle\\ \langle z,y\rangle & \langle z,z \rangle & \langle z,w\rangle & \langle z,v\rangle\\ \langle w,y \rangle & \langle w,z \rangle & \langle w,w\rangle & \langle w,v\rangle\\\langle v,y \rangle & \langle v,z \rangle & \langle v,w\rangle & \langle v,v\rangle
\end{vmatrix}\nonumber\\
&=&\langle x,y|v,w,z\rangle\Vert w|z\Vert^2\Vert z\Vert^4.\label{ee5}
\end{eqnarray}

From relation (\ref{ee5}), for $x=y$ we deduce:
\begin{equation}\label{ee6}
(x,x|v,w,z)_*=\Gamma(x,v,w,z)\cdot\Vert w|z\Vert^{2}\Vert z\Vert^4,
\end{equation}
where $\Gamma(x,v,w,z)$ is the Gram's determinant. 
\par In \cite{2}, Cho, Mati\'c and Pe\u cari\'c, used Gram's determinant $\Gamma(x_1,x_2,...,x_k|z)$ of the vectors $x_1, x_2,..., x_k$ with respect to the vector $z$ by:
\begin{equation}\label{ee7}
\Gamma(x_1,x_2,...,x_k|z)= \begin{vmatrix}
\langle x_1,x_1|z\rangle & \langle x_1,x_2 |z\rangle &...& \langle  x_1,x_k|z\rangle\\ \langle  x_2,x_1|z\rangle & \langle  x_2,x_2 |z\rangle &...& \langle  x_2,x_k|z\rangle \\ \vdots &\vdots &\ddots &\vdots\\ \langle  x_k,x_1 |z\rangle & \langle   x_k,x_2 |z\rangle &...& \langle   x_k,x_k|z\rangle
\end{vmatrix}.
\end{equation}

We consider the following determinant, which can be rewritten using formula (\ref{e41a}):
\begin{eqnarray}
&&\begin{vmatrix}
\langle x,y|z\rangle & \langle  x,w |z\rangle & \langle  x,v|z\rangle \\ \langle  w,y|z\rangle & \langle  w,w |z\rangle & \langle  w,v|z\rangle \\ \langle  v,y |z\rangle & \langle  v,w |z\rangle &\langle v,v|z\rangle
\end{vmatrix}\nonumber\\
&=&\frac{1}{\langle w,w|z\rangle}\Big[(x,y|w,z)_*(v,v|w,z)_*-(x,v|w,z)_*(v,y|w,z)_*\Big]\nonumber\\ 
&=&\frac{1}{\langle w,w|z\rangle}(x,y|v,w,z)_*.\label{ee9}
\end{eqnarray}

From relations (\ref{ee5}) and (\ref{ee9}) we find the following identity:\\
\begin{equation}
\begin{vmatrix}
\langle x,y|z\rangle & \langle  x,w |z\rangle & \langle  x,v|z\rangle \\ \langle  w,y|z\rangle & \langle  w,w |z\rangle & \langle  w,v|z\rangle \\ \langle  v,y |z\rangle & \langle  v,w |z\rangle & \langle  v,v|z\rangle
\end{vmatrix}=\langle x,y|v,w,z\rangle\Vert z\Vert^4
\end{equation}
which implies the relation:

\begin{equation}
\Gamma(x,w,v|z)=\Gamma(x,w,v)\Vert z\Vert^4.
\end{equation}
\par\

{\bf 4.} Let $x, y, e,w$ be vectors in the inner product space $X,$ over the field of real numbers and the vectors $\{e,x,y\}$ being linearly independent, such that $$ax+by+ce=w,$$ where $a,b,c\in\mathbb{R}.$ 

We want to study the problem of determining the scalars $a, b, c$.
Using the inner product and its properties, we deduce

\begin{equation}
\left\{\begin{matrix} 
& a\langle x,x\rangle+b\langle y,x\rangle+c\langle e,x\rangle=\langle w,x\rangle\\
& a\langle x,y\rangle+b\langle y,y\rangle+c\langle e,y\rangle=\langle w,y\rangle\\
& a\langle x,e\rangle+b\langle y,e\rangle+c\langle e,e\rangle=\langle w,e\rangle.\end{matrix}\right.
\end{equation}
Therefore, we have to solve this system with three equations and three unknowns $a,b,c\in\mathbb{R}.$ The matrix of the system is
$$A=\begin{pmatrix} \langle x,x \rangle \ \langle y,x \rangle \ \langle e,x \rangle \\ \langle x,y \rangle \ \langle y,y \rangle \ \langle e,y \rangle\\ \langle x,e \rangle \ \langle y,e \rangle \  \langle e,e \rangle
\end{pmatrix}.$$

From formula (\ref{ee3}) we find:

$$\det A=\Gamma(x,y,e)=\frac1{\Vert e\Vert^2}(x,x|y,e)_*.$$

\par From P1) $(x,x|y,e)_*$  is zero if and only if the vectors $x,y,e$ are linearly dependent. But, the vectors $\{e,x,y\}$ are linearly independent, therefore, we have $(x,x|y,e)_*>0$.
Using the Cramer method, we find that
$$
a=\frac{( w,x|y,e)_*}{(x,x|y,e)_*}, \
b=\frac{(w,y|x,e)_*}{(x,x|y,e)_*}, \
c=\frac{\Vert e\Vert^2(w,e|x,y)_*}{\Vert y\Vert^2(x,x|y,e)_*}. 
$$

\par\ 

In the particular case when $\Vert e\Vert=1$ we obtain:
$$ a=\frac{( w,x|y,e)_*}{(x,x|y,e)_*}, \ b=\frac{(w,y|x,e)_*}{(x,x|y,e)_*}, \ c=\langle w,e\rangle-a\langle x,e\rangle-b\langle y,e\rangle.$$

\par {\bf 5.} Next, we will make a correlation of the previous calculations with the coefficients that appear in the case of a multiple linear regression model.\\
\indent A process is called \textit{multiple linear regression}, when we have more than one independent variable \cite{7}. For a general linear model for two independent variables $V$ and $W$ and a dependent variable $Z$, $Z=aV+bW+c$, where    
$V=\begin{pmatrix}
x_{i}\\
\frac{1}{n}
\end{pmatrix}_{1\leq i\leq n}$; 
$W=\begin{pmatrix}
y_{i}\\
\frac{1}{n}
\end{pmatrix}_{1\leq i\leq n}$; 
$Z=\begin{pmatrix}
z_{i}\\
\frac{1}{n}
\end{pmatrix}_{1\leq i\leq n}$
with probabilities $P(V=x_{i})=\dfrac{1}{n}$, $P(W=y_{i})=\dfrac{1}{n} $, $P(Z=z_{i})=\dfrac{1}{n},$ for any $i=\overline{1,n}.$  \\ 
\indent We can describe the underlying relationship between $z_i$ and $x_i, y_i$ involving error term $\epsilon_i$ by  $\epsilon_i=z_i-ax_i-by_i-c$.\\ 
\indent If we take $S(a,b,c)=\sum_{i=1}^{n}\epsilon_{i}^2=\sum_{i=1}^{n}(z_i-ax_i-by_i-c)^2$, then we have to find $\underset{a,b,c\in\mathbb{R}}{min}S(a,b,c).$   Using the Lagrange method, we obtain

\begin{eqnarray*} 
&&a\sum_{i=1}^{n}x_i+b\sum_{i=1}^{n}y_i+nc=\sum_{i=1}^{n}z_i,\\ &&a\sum_{i=1}^{n}x_i^2+b\sum_{i=1}^{n}x_iy_i+c\sum_{i=1}^{n}x_i=\sum_{i=1}^{n}x_iz_i,\\ &&a\sum_{i=1}^{n}x_iy_i+b\sum_{i=1}^{n}y_i^2+c\sum_{i=1}^{n}y_i=\sum_{i=1}^{n}y_iz_i. 
\end{eqnarray*}
By simple calculations, we deduce 

\begin{eqnarray*}
a&=&\frac{Var(W)Cov\left( V,Z\right)-Cov(V,W)Cov(W,Z)}{Var\left( V\right)Var(W)-Cov^2(V,W)},\\
b&=&\frac{Var(V)Cov\left( W,Z\right)-Cov(V,W)Cov(V,Z)}{Var\left( V\right)Var(W)-Cov^2(V,W)},\\
c&=&E\left(Z\right)-aE\left(V\right)-bE\left(W\right).
\end{eqnarray*}
Now, we take the vector space $(X=\mathbb{R}^n,\langle\cdot,\cdot\rangle).$ For $x=(x_{1},x_{2},...,x_{n}), \ y=(y_{1},y_{2},...,y_{n}), \ z=(z_{1},z_{2},...,z_{n}),$ we have 
$$\langle x,y \rangle=x_{1}y_{1}+x_{2}y_{2}+...+x_{n}y_{n}, \ \Vert x\Vert=\sqrt{x_{1}^2+x_{2}^2+...+x_{n}^2},$$

$$(x,y|z)_*=\langle x,y \rangle \langle z,z \rangle - \langle x,z \rangle\langle z,y \rangle=\sum_{i=1}^{n}x_{i}y_{i}\sum_{i=1}^{n}z_{i}^2-\sum_{i=1}^{n}x_{i}z_{i}\sum_{i=1}^{n}z_{i}y_{i}$$  and
$$ \Vert x|z \Vert =\sqrt{\sum_{i=1}^{n}x_{i}^2\sum_{i=1}^{n}z_{i}^2-\bigg(\sum_{i=1}^{n}x_{i}z_{i}\bigg)^2}.$$
If $e=\frac{u}{\Vert u\Vert}$, where $u=(1,1,...,1)\in\mathbb{R}^n$, then the \textit{average} of vector $x$ is $\mu_{x}=\bigg\langle\frac{x}{\Vert u\Vert},e\bigg\rangle=\frac{1}{n}\sum_{i=1}^{n}x_{i}$, and we have
$$\bigg \Vert \frac{x}{\Vert u\Vert}|e\bigg\Vert=\sqrt{\frac{1}{n}\sum_{i=1}^{n}x_{i}^2-\bigg( \frac{1}{n}\sum_{i=1}^{n}x_{i}\bigg)^2}.$$

Therefore, in $(\mathbb{R}^n,\langle \cdot,\cdot \rangle)$, we define the \textit{variance} of a vector $x$ by $var(x):=\bigg\Vert\frac{x}{\Vert u\Vert}|e\bigg\Vert^2.$

The \textit{standard deviation} $\sigma(x)$ of $x\in \mathbb{R}^n$ is defined by 
$\sigma(x):=\sqrt{var(x)}$, so we deduce that $\sigma(x)=\bigg\Vert\frac{x}{\Vert u\Vert} |e\bigg\Vert.$ 
Since, using the standard $2$-inner product, we have
$$\bigg(\frac{x}{\Vert u\Vert},\frac{y}{\Vert u\Vert}|e\bigg)_*=\frac{1}{n}\sum_{i=1}^{n}x_{i}y_{i}-\bigg(\frac{1}{n}\sum_{i=1}^{n}x_{i}\bigg)\bigg(\frac{1}{n}\sum_{i=1}^{n}y_{i}\bigg),$$
it is easy to define the \textit{covariance} of two vectors $x$ and $y$ by 
 $$cov(x,y):=\bigg(\frac{x}{\Vert u\Vert},\frac{y}{\Vert u\Vert}|e\bigg)_*.$$
It is easy to see that, we obtain 

\begin{eqnarray*}
a&=&\frac{var(y)cov\left( x,z\right)-cov(x,y)cov(y,z)}{var\left( x\right)var(y)-cov^2(x,y)},\\
b&=&\frac{var(x)cov\left( y,z\right)-cov(y,x)cov(x,z)}{var\left( x\right)var(y)-cov^2(x,y)},\\
c&=&\mu_z -a\mu_x-b\mu_y.
\end{eqnarray*}
We observe that by the vector method we obtain the same coefficients as by the Lagrange method.\\

\par {\bf 6.} In \cite{13}, the \textit{Chebyshev functional} is defined by 
$$T_{z}(x,y)=\Vert z\Vert^2\langle x,y\rangle-\langle x,z\rangle\langle y,z\rangle,$$
for all $x,y\in X$, where $z\in X$ is a given nonzero vector.

It is easy to see that if the standard $2$-inner product $(\cdot,\cdot |\cdot)$ is defined by the inner product $\langle\cdot,\cdot\rangle,$ then we have  $T_{z}(x,y)=(x,y|z)_*=(x,y|z)$.

Therefore, we generalize this Chebyshev functional to the following functional:
$$T_{x_{n},...,x_{2}}(x,y):=(x,y|x_{n},...,x_{2})_*,$$
which we will call \textit{n-Chebyshev functional}, so

\begin{eqnarray}
T_{x_{n},...,x_{2}}(x,y)&=&T_{x_{n-1},...,x_{2}}(x,y)T_{x_{n-1},...,x_{2}}(x_{n},x_{n})\nonumber\\
&&\qquad-T_{x_{n-1},...,x_{2}}(x,x_{n})T_{x_{n-1},...,x_{2}}(x_{n},y),\label{ee11}
\end{eqnarray}
for all $x,y\in X$, where $x_{2},...,x_{n}\in X$  are given nonzero vectors.

\par In a particular case, when $n=3$, we have
$$ T_{w,z}(x,y)=(x,y|w,z)_*=(x,y|z)_*(w,w|z_*)-(x,w|z)_*(w,y|z)_*$$ 
so, we have
\begin{eqnarray*}
T_{w,z}(x,x)&=&(x,x|w,z)_*=(x,x|z)_*(w,w|z)_*-(x,w|z)_*(w,x|z)_*\\ 
&=&\Vert x|z\Vert^2\Vert w|z\Vert^2-(x,w|z)^2\\
&=&(\Vert x\Vert^2\Vert w\Vert^2\Vert z\Vert^2+2\langle w,z\rangle\langle z,x\rangle\langle x,w\rangle
 -\Vert x\Vert^2\langle w,z\rangle^2-\Vert w\Vert^2\langle z,x\rangle^2\\
&&-\Vert z\Vert^2\langle x,w\rangle^2)\Vert z\Vert^2.
\end{eqnarray*}

Therefore, the Cauchy-Schwarz inequality in terms of the $n$-Chebyshev functional becomes: 
\begin{equation}
|T_{x_{n},...,x_{2}}(x,y)|^2\leq T_{x_{n},...,x_{2}}(x,x)T_{x_{n},...,x_{2}}(y,y). 
\end{equation}          

\section{Conclusions}\label{sec:5}
\indent In this paper we exemplified the weak $n$-inner product only by the weak $n$ iterated $2$-inner product. This particular case of  weak $n$-inner product  does not exhaust all the possibilities of particular cases. The weak $n$-inner product is clearly more general then the $n$-inner product and consequently it offers more possibilities. An important connection is between the vector method and the Lagrange method given above. In the future, we will determine a formula for multiple regression for $n$ independent variables.

\textbf{Acknowledgment.} 
The authors would like to thank to the reviewers for the pertinent remarks, which led to an improvement of the paper.

\par\ \\
Department of Mathematics and Computer Science, \\
{\it Transilvania} University of Bra\c{s}ov \\
Iuliu Maniu 50, Bra\c sov, Romania,\\
E-mail: minculeten@yahoo.com
\par\ \\
Department of Mathematics and Computer Science, \\
{\it Transilvania} University of Bra\c{s}ov \\
Iuliu Maniu 50, Bra\c sov, Romania,\\
  E-mail: radupaltanea@yahoo.com

\end{document}